\documentclass[12pt]{amsart}
\usepackage{amscd,amssymb}
\usepackage[graph,frame,poly,arc]{xy}  
\usepackage[plainpages,backref,urlcolor=blue]{hyperref}

\theoremstyle{plain}
\newtheorem{thm}[subsection]{Theorem}
\newtheorem{lem}[subsection]{Lemma}
\newtheorem{prop}[subsection]{Proposition}
\newtheorem{cor}[subsection]{Corollary}

\theoremstyle{definition}
\newtheorem{rk}[subsection]{Remark}
\newtheorem{definition}[subsection]{Definition}
\newtheorem{ex}[subsection]{Example}
\newtheorem{conj}[subsection]{Conjecture}

\numberwithin{equation}{section}
\newcommand{\OO}{{\mathcal O}}

\newcommand{\A}{{\mathcal A}}
\newcommand{\B}{{\mathcal B}}

\newcommand{\Z}{\mathbb{Z}}

\newcommand{\C}{\mathbb{C}}
\newcommand{\K}{\mathbb{K}}
\newcommand{\PP}{\mathbb{P}}

\DeclareMathOperator{\coker}{coker}

\DeclareMathOperator{\mult}{mult}

\DeclareMathOperator{\Der}{Der}

\DeclareMathOperator{\POexp}{POexp}

\begin{document}

\title [Addition-deletion results for logarithmic derivations]
{Addition-deletion results for the minimal degree of logarithmic derivations of arrangements}

\author[Takuro Abe]{Takuro Abe}
\address{Institute of Mathematics for Industry,
Kyushu University,
Fukuoka 819-0395, Japan.}
\email{abe@imi.kyushu-u.ac.jp}

\author[Alexandru Dimca]{Alexandru Dimca}
\address{Universit\'e C\^ ote d'Azur, CNRS, LJAD, France and Simion Stoilow Institute of Mathematics,
P.O. Box 1-764, RO-014700 Bucharest, Romania}
\email{dimca@unice.fr}

\author[Gabriel Sticlaru]{Gabriel Sticlaru}
\address{Faculty of Mathematics and Informatics,
Ovidius University
Bd. Mamaia 124, 900527 Constanta,
Romania}
\email{gabrielsticlaru@yahoo.com }




\begin{abstract} 
 We study the change of the minimal degree of a logarithmic derivation  of a hyperplane arrangement  under the addition or the deletion of a hyperplane, and give a number of applications. First, we prove the existence of Tjurina maximal line arrangements in a lot of new situations.
Then, starting with Ziegler's example of a pair of arrangements of $d=9$ lines with $n_3=6$ triple points in addition to some double points, having the same combinatorics, but distinct minimal degree of a logarithmic derivation, we construct new examples of such pairs, for any number $d\geq 9$ of lines, and any number $n_3\geq 6$ of triple points. Moreover, we show that such examples are not possible for line arrangements having only double and triple points, with $n_3 \leq 5$. 

\end{abstract}
 
\maketitle


\section{Introduction} 
Let $\K$ be a field of characteristic zero, consider the polynomial ring
$S=\K[x_1,\ldots,x_{\ell}]$ with the usual grading, i.e., 
$S=\oplus_{d \in \Z_{\ge 0}} S_d$, and for an $S$-graded module $M$, 
let $M=\oplus_{d \in \Z} M_d$ be its decomposition according to the grading.
Let $X$ be a reduced projective hypersurface
in $\PP^{\ell -1}$, defined by a homogeneous polynomial
$f\in S_d$ of degree $d$. We assume that $X$ is \textbf{essential}, that is $X$ is not the cone over a projective hypersurface in some $\PP^n$ with $n <\ell-1$. When $X=\A$ is a \textbf{hyperplane arrangement}, the main situation considered below, this definition agrees with the usual one. The details will be recalled in the next section.
Let 
$$
\Der S:=\oplus_{i=1}^\ell S \partial_{x_i}
$$ be  the module of 
derivations of $S$,  a $\Z$-graded free $S$-module of rank $\ell$. Here $0 \neq \theta \in \Der S$ is homogeneous of degree 
$e$ if $\theta(g)$ is zero or homogeneous of degree $e$ for all $g \in S_1$. 
For example, the  \textbf{Euler derivation} 
$\theta_E:=\sum_{i=1}^\ell x_i \partial_{x_i}$ is homogeneous of degree $1$. 
The \textbf{logarithmic derivation module} $D(X)$ 
of the hypersurface $X$ is defined by 
$$
D(X):=\{\theta \in \Der S \mid 
\theta(f) \subset (f)\},
$$
where $(f)$ denotes the principal ideal generated by $f$ in $S$.
It is known that $D(X)$ is an $S$-graded reflexive module, but not free in general. It is clear that $\theta_E \in D(X)$.
We say that 
$X$ has  \textbf{exponents} $\exp(X)=(d_1,\ldots,d_k)$ if 
there are homogeneous derivations $\theta_1=\theta_E,\ldots,\theta_k$ with
$\deg \theta_j=d_j$ which form a minimal set of generators for the graded $S$-module $D(X)$. Since $X$ is essential, it follows that $d_j>0$ for all $j$.
When these integers $d_j$ are written in increasing order, we use the notation
$$(d_1,\ldots,d_k)_{\leq}.$$
 Consider the graded $S$-submodule
$$
D_0(X)=\{
\theta \in D(X) \mid \theta (f) =0\}.
$$
in $D(X)$ and note the decomposition
$$
D(X)=D_0(X) \oplus S \theta_E.$$
Because of this decomposition, it is usual to choose the minimal generators $\theta _j$ above such that $\theta_j \in D_0(X)$ for $j>1$.
If $\theta=\sum_{i=1}^{\ell}a_i\partial_{x_i}$ with $a_i \in S_r$ for some integer $r$, the condition $\theta (f) =0$ translates into the following homogeneous \textbf{Jacobian relation} or  \textbf{Jacobian syzygy}
\begin{equation}
\label{syzy}
\sum_{i=1}^{\ell}a_if_{x_i}=0,
\end{equation}
involving the partial derivatives $f_{x_i}=\partial_{x_i}f$ of the polynomial $f$. In this way, the generators $\theta_j$ for $j>1$ are sometimes identified with Jacobian relations. This explains the following.
\begin{definition}
The \textbf{minimal degree of a Jacobian relation} of $X$, denoted by $r(X)$ or $mdr(X)$, is defined by 
$$r(X)=\min_{r \in \Z} \{r \mid D_0(X)_r \neq (0)\}.$$
\label{defmdrX}
\end{definition}
If $(d_1,\ldots,d_k)_{\leq}$ are the exponents of $X$, then $d_1=1$ and
$d_2=r(X)$. We say that the hypersurface $X$ is \textbf{free} if $k=\ell$,
i.e. the graded $S$-module $D(X)$ is free. When this happens, one has
$$d_1+\ldots+d_{\ell}=d.$$
In particular, for a free plane curve $X$ (the case when $\ell=3$), 
the exponents are determined by $r(X)$, namely
\begin{equation}
\label{expC}
\exp(X)=(1,r(X),d-1-r(X)),
\end{equation}
and $r(X) \leq (d-1)/2$.  Recall that a plane curve $X$ is \textbf{nearly free} when its exponents are given by 
$\exp(X)=(1,r(X),d-r(X),d-r(X))$ with the unique relation 
at degree $d-r(X)+1$, see \cite{DStExpo,DSt3syz}.

The main motivation of this paper, and the reason to study the invariant $r(X)$, is the following conjecture due to H. Terao.
\begin{conj}
\label{conjT1}
Let $\A$ and $\B$ be two hyperplane arrangements, having isomorphic intersection lattices $L(\A) \cong L(\B)$. If $\A$ is free, then $\B$ is also free.
\end{conj}
For more on Terao's conjecture, as well as for basic information on hyperplane arrangements, we refer to \cite{DHA, OT}.
This conjecture is open, even in the case of line arrangements in $\PP^2$, in spite of a lot of work and partial results in the recent years, see \cite{A2,A4,DIM,Y3}.
Note that the freeness of a line arrangement $\A$ is \textbf{not determined} by the \textbf{weak combinatorics} of $\A$, namely the numbers $n_j$ of points in $\A$ of multiplicity $j\geq 2$, see \cite{MV19}.
In the case of line arrangements, using \eqref{expC} and a result by A. du Plessis and C.T.C. Wall quoted below in Theorem \ref{thmCTC}, Terao's conjecture can be restated as follows.
\begin{conj}
\label{conjT2}
Let $\A$ and $\B$ be two line arrangements, having isomorphic intersection lattices $L(\A) \cong L(\B)$. If $\A$ is free, then $r(\A)=r(\B)$.
\end{conj}
It is known that the intersection lattice $L(\A)$ does not determine the integer $r(\A)$ in general: indeed, G. Ziegler produced two arrangements $\A$ and $\B$ of $d=9$ lines, having only double and triple points, such that $L(\A) \cong L(\B)$, and $5=r(\A) \ne r(\B)=6$, see Remark \ref{rk1.5} for more details. However, the following \textbf{stronger form of Terao's conjecture} might be true.
\begin{conj}
\label{conjT3}
Let $\A$ and $\B$ be two arrangements of $d$ lines, having isomorphic intersection lattices $L(\A) \cong L(\B)$. If $r(\A)<d/2$, then $r(\A)=r(\B)$.
\end{conj}
Note that in \cite{MV19}, the authors produce two arrangements $\A$ and $\B$ of $d$ lines, having the same weak combinatorics, and such that
$r(\A) <d/2$ and $r(\A) \ne r(\B)$.
Conjecture \ref{conjT3} can be stated in a more geometric way as follows, when $\K=\C$. Let $E(\A)$ be the rank 2 vector bundle on $\PP^2$ naturally associated with the reflexive graded $S$-module $D_0(\A)$.
For a generic line $L$ in $\PP^2$, the restriction $E(\A)|L$ splits as a direct sum $\OO_L(-e_1)\oplus\OO_L(-e_2)$. The pair $(e_1,e_2)$ is called the \textbf{generic splitting type of the bundle} $E(\A)$ and it is known that the two pairs $(e_1,e_2)$ and $(d,r(\A))$ determine each other when $r(\A)<d/2$, with $d=|\A|$, see \cite[Propositions 3.1 and 3.2]{AD} and \cite[Theorem 1.2]{Drcc}. When $r(\A)\geq d/2$, it follows from \cite{AD, Drcc} that the generic splitting type $(e_1,e_2)$ is determined by $d$ and the global Tjurina number $\tau(\A)$, which is determined in turn by the 
weak combinatorics of $\A$ via the well known formula
\begin{equation}
\label{TauW}
\tau(A)=\sum_{j\geq 2} n_j(j-1)^2.
\end{equation}
Hence Conjecture \ref{conjT3} is equivalent, when $\K=\C$, to the following conjecture, which has already appeared in \cite[Question 7.12]{CHMN} and in \cite{AD}.
\begin{conj}
\label{conjT4}
Let $\A$ and $\B$ be two line arrangements, having isomorphic intersection lattices $L(\A) \cong L(\B)$. 
Then the rank 2 vector bundles $E(\A)$ and $E(\B)$ have the same
generic splitting type.
\end{conj}
Note that in spite of Ziegler's example mentioned above, Conjecture \ref{conjT4} holds for line arrangements having only double and triple points, see Remark \ref{rkNU}.

In this paper we start a detailed investigation of the dependence of the minimal degree $r(\A)$ of a Jacobian relation  of a hyperplane arrangement $\A$ on the combinatorics of $\A$. As a first step, we
study the change of the invariant $r(\A)$ of a hyperplane arrangement $\A$ under the addition or the deletion of a hyperplane $H$, and give a number of applications.

In section 2, after some preliminary material on arrangements, we establish the main general {\bf addition-deletion result for the invariant} $r(\A)$ of a hyperplane arrangement $\A$, see Theorem \ref{ADgeneral}. The special case of free hyperplane arrangements is discussed in Theorem \ref{freemdr}. Other authors have considered addition-deletion to study the logarithmic derivation module $D(\A)$, see for instance \cite{STY,W},
but without paying attention to the invariant $r(\A)$.

In section 3 we concentrate our attention to a line arrangement $\A$ in $\PP^2$. The corresponding addition-deletion results for $r(\A)$ are stated in Theorems \ref{main} and \ref{main2}, while the case of free line arrangements is discussed in Theorem \ref{freemdr3}. We then recall the relation between the invariant $r(\A)$ and the maximal multiplicity $m(\A)$ of an intersection point of the line arrangement $\A$ following \cite{DMich}. Corollary \ref{cor1} says that $r(\A)$ is determined by the weak combinatorics of $\A$ when $2m(\A) \geq |\A|$.

In section 4 we give some applications of the above results. The result by A. du Plessis and C.T.C. Wall quoted below in Theorem \ref{thmCTC} gives an upper bound $\tau(d,r)_{max}$ for the global Tjurina number $\tau(C)$ of a reduced plane curve $C$, in terms of its degree $d$ and the invariant $r=r(C)$. A curve  $C$, for which the equality $\tau(C)=\tau(d,r)_{max}$ holds, is called a \textbf{maximal Tjurina curve of type} $(d,r)$.
For any pair
$(d,r)$, with $1 \leq r <d/2$, a maximal Tjurina curve of type $(d,r)$ is nothing else but a free curve $C$ of degree $d$  with $r(C)=r$, and the existence of such curves, even in the class of line arrangements, follows from \cite{DStExpo}. 
For the pairs of the form $(d=2r,r)$, a maximal Tjurina curve of type $(d,r)$ is nothing else but a nearly free curve $C$ of degree $d=2r$  with $r(C)=r$, and the existence of such curves, even in the class of line arrangements, follows again from \cite{DStExpo}. The existence of maximal Tjurina curves of type $(d,r)$ when $d/2<r\leq d-2$, is much more subtle. The following conjecture was stated in \cite{DS}.
\begin{conj}
\label{conjMT1}  For any integer $d \geq 5$ and for any integer $r$ such that $d/2 < r\leq d-2$, there are maximal Tjurina line arrangements of type $(d,r)$.
\end{conj}
As noted in \cite{DS}, the generic line arrangement of $d$ lines is Tjurina maximal of type $(d,d-2)$ for any $d \geq 2$, see also Remark \ref{rkSE} below.
Line arrangements which are {\bf potentially}  maximal Tjurina  of the following types:
\begin{itemize}
\item[(1)] $(d,r)=(2r-1,r)$ for $r \geq 3$,

\item[(2)] $(d,r)=(d,d-4)$ for $d \geq 8$, and

\item[(3)] $(d,r)=(d,d-3)$ for $d \geq 7$

\end{itemize}
have been put forth in \cite{DS}, following numerical experiments with SINGULAR. The fact that these arrangements are indeed maximal Tjurina is proved here, see Corollary \ref{corex5.5} for type (1), Theorem \ref{d-4} for type (2), and Theorem \ref{d-3} for type (3). As a result, Conjecture \ref{conjMT1} holds in all these extremal cases for $r$ in the interval $d/2<r< d-2$. The existence of maximal Tjurina curves of type $(d,r)$ when $d/2<r< d-2$, in a lot of new cases, is proved in 
Theorem \ref{thm100}, Proposition \ref{prop101} and Remark \ref{rk101}.

We continue section 4 by investigating the effect on $r(\A)$ of adding a generic line, either passing through a point of maximal multiplicity of $\A$, or just transversal to $\A$, see Propositions \ref{prop2} and \ref{prop4}. Using these results and our main addition-deletion result, Theorem \ref{main},
we determine the invariant $r(\A)$ in the case of line arrangements having only double and triple points, when the number of triple points $n_3$ is $\leq 5$, see Theorem \ref{thm2}. The conclusion is that in these cases
the invariant $r(\A)$ is determined by the combinatorics of $\A$, in a precise, but rather complicated way. Ziegler's example, which was the only known example of this type until now, shows that this result is sharp, i.e. it does not extend for the situation $n_3 \geq 6$, see Corollary \ref{cor2}. In fact, using Ziegler's example and adding well chosen lines, we can construct similar examples of pairs of arrangements of $d$ lines, having only double and triple points, for any $d \geq 9$ and any possible weak combinatorial data $(n_2,n_3)$, when $n_3 \geq 6$, see the proof of Corollary \ref{cor2}.
\medskip

\noindent
\textbf{Acknowledgements}. The first author is partially supported by 
KAKENHI, 
Fund for the Promotion of Joint International Research (Fostering Joint International Research (A)) 18KK0389.
The second author has been supported by the French government, through the $\rm UCA^{\rm JEDI}$ Investments in the Future project managed by the National Research Agency (ANR) with the reference number ANR-15-IDEX-01 and by the Romanian Ministry of Research and Innovation, CNCS - UEFISCDI, grant PN-III-P4-ID-PCE-2016-0030, within PNCDI III.

\section{Hyperplane arrangements} 

\subsection{Preliminaries}
 First we recall some definitions and notations.
Let $V=\K^\ell$, $x_1,\ldots,x_\ell$ a basis for $V^*$ and let 
$S:=\mbox{Sym}^*(V^*)=\K[x_1,\ldots,x_\ell]$. 
We say that $\A$ is a \textbf{hyperplane arrangement} in $V$ if $\A$ is a finite set of linear hyperplanes in $V$. We say that $\A$ is \textbf{essential} if $\cap_{H \in \A} H 
=\{0\}$. 
We assume that all arrangements are essential unless otherwise specified. 
For $H \in \A$, let 
$$
\A^H:=\{H \cap L \mid L \in \A \setminus \{H\}\}
$$
be the \textbf{restriction}.
Let 
$$
L(\A):=\{
\cap_{H \in \B} H \mid 
\B \subset \A\}
$$ be the \textbf{intersection lattice} of $\A$. Then we can define the \textbf{M\"{o}bius function} $\mu:L(\A) \rightarrow \Z$ by 
$\mu(V)=1$, and by 
$$
\mu(X):=-\sum_{X \subsetneq Y \subset V,\ Y \in L(\A)} \mu(Y).
$$
Then we can define the \textbf{characteristic polynomial} $\chi(\A;t)$ by 
$$
\chi(\A;t):=\sum_{X \in L(\A)} \mu(X)t^{\dim X}=\sum_{i=0}^\ell 
b_i(\A)t^i.
$$
If $\A \neq \emptyset$, then $\chi(\A;t)$ is divisible by $(t-1)$. Let 
$$
\chi_0(\A;t):=\chi(\A;t)/(t-1)=\sum_{i=0}^{\ell-1} b_i^0(\A) t^i.
$$
It is easy to show that 
$$
b_1(\A)=|\A|, b_1^0(\A)=|\A|-1,\ b_2^0(\A)=b_2(\A)-|\A|+1.
$$
For $H \in \A$ fix a linear form $\alpha_H \in V^*$ such that 
$\ker \alpha_H=H$. Then the \textbf{logarithmic derivation module} $D(\A)$ can be defined in this situation as follows:
$$
D(\A):=\{\theta \in \Der S \mid 
\theta (\alpha_H) \in S\alpha_H\ (\forall H \in \A)\}.
$$
For $Q(\A):=
\prod_{H \in \A} \alpha_H$, one has as in the Introduction
$$
D_0(\A):=\{
\theta \in \Der S \mid \theta (Q(\A))=0\}.
$$
The first easy, but important lemma is the following.

\begin{lem}[Lemma 1.33, \cite{Y3}]
For $H \in \A$, let $D_H(\A):=\{ 
\theta \in D(\A) \mid \theta(\alpha_H)=0\}$. Then 
$$
D(\A)=D_0(\A) \oplus S \theta_E=
D_H(\A) \oplus S \theta_E.$$
In particular, if $\A 
\neq \emptyset$, 
$$
D_0(\A) \simeq D(\A)/S \theta_E \simeq D_H(\A)
$$
for any $H \in \A$.
\label{lemma1}
\end{lem}

Lemma \ref{lemma1} is a well-known classical result in arrangement theory. It implies in particular the equality
\begin{equation}
r(\A)=\min_{d \in \Z} \{d \mid D_H(\A)_d \neq (0)\},
\label{defmdr}
\end{equation}
for any $H \in \A$. 
In the study of $r(\A)$, Lemma \ref{lemma1} shows a big difference of hyperplane arrangements compared with, say,  the case of general plane curves. The reason is that, for $L\in \A':=\A \setminus \{H\}$, Lemma \ref{lemma1} and the definition of 
logarithmic vector fields show that 
$$
D_0(\A) \simeq D_L(\A) \subset D_L(\A') \simeq D_0(\A').
$$
Thus we can directly compare $r(\A)$ and 
$r(\A')$. To compare them more precisely, the following result due to Terao always plays the key role.

\begin{thm}[Terao's polynomial $B$-theory, \cite{T1}]
Let $H \in \A$, 
$\A':=\A 
\setminus \{H\}$, 
and let us define the homogeneous degree 
$|\A'|-|\A^H|$-polynomial $B$ by 
$$
B:=\prod_{X \in \A^H} \alpha_{\nu(X)},
$$
where $\nu:\A^H \rightarrow \A$ is a section 
satisfying that 
$\nu(X) \cap H=X$. Then 

(1)\,\,
for an arbitrary $\theta \in D(\A')$, it holds that 
$$
\theta(\alpha_H) \in (\alpha_H,B),
$$
where $(\alpha_H,B)$ denotes the ideal of $S$ generated by $\alpha_H$ and $B$. 
Thus, $\theta \in D(\A')$ is in $D(\A)$ if $\deg \theta <|\A'|-|\A^H|$.

(2)\,\ Assume that there is $\varphi \in 
D(\A')$ such that $\deg \varphi=|\A'|-|\A^H|$ and that 
$\varphi \not \in D(\A)$. Then for $\theta \in D(\A')$, there is 
$f \in S$ such that $\theta -f \varphi \in D(\A)$. Thus 
$$
D(\A')=D(\A)+S\cdot \varphi.
$$
\label{B}
\end{thm}

To compare algebraic structures of $D(\A)$ and $D(\A')$, the most useful tool is 
Terao's addition-deletion. Since $r(\A)$ sees only the lowest degree generator of $D_0(\A)$, 
the following variant of the addition-deletion theorem is useful.

\begin{thm}[Multiple deletion theorem, \cite{AT}]
Let $\A$ be a free hyperplane arrangement with $\exp(\A)=(1,d_2,\ldots,d_\ell)_\le$.  If there is 
$H \in \A$ such that $|\A|-|\A^H|=d_2$, then 
$\A':=\A \setminus 
\{H\}$ is free with $\exp(\A')=(1,d_2-1,d_3,\ldots,d_\ell)_\le$. 
\label{MDT}
\end{thm}

To compare $r(\A)$, the following two restriction maps play important roles. Let us introduce them. First, the \textbf{Euler restriction} $\rho:D(\A) \rightarrow D(\A^H)$ 
is defined by 
$$
\rho(\theta)(\overline{f}):=\overline{\theta(f)}
$$
for $\theta \in D(\A),\ f \in S/\alpha_H S$. Here $\overline{f}$ denotes the 
image of $f \in S$ by the canonical surjection $S \rightarrow 
S/\alpha_H S$. It is well-known that there is an exact sequence 
$$
0 \rightarrow D(\A \setminus \{H\}) \stackrel{\cdot \alpha_H }{\rightarrow }
D(\A) \stackrel{\rho}{\rightarrow} D(\A^H).
$$

Also, we have the other restriction. To introduce it, let us recall 
multiarrangements.
For an arrangement $\A$, let $m:\A \rightarrow \Z_{>0}$ be a \textbf{multiplicity}, and the pair 
$(\A,m)$ is called a \textbf{multiarrangement}. 
For $H \in \A$, let $\delta_H$ be a multiplicity on $\A$ defined by $\delta_H(L)=1$ if $L=H$, and $0$ otherwise. 
We can define its \textbf{logarithmic derivation module} $D(\A,m)$ by 
$$
D(\A,m):=\{\theta \in \Der S \mid 
\theta(\alpha_H) \in S \alpha_H^{m(H)}\ (\forall H \in \A)\}.
$$
We can define the \textbf{freeness} and 
\textbf{exponents} of $D(\A,m)$ in the same manner as for $\A$. We can construct multiarrangements canonically from an 
arrangement $\A$ and $H \in \A$ as follows. Define the \textbf{Ziegler multiplicity} $m^H$ on $\A^H$ by 
$$
m^H(X):=|\{
L \in \A \setminus \{H\} \mid 
L \cap H=X\}
$$
for $X \in \A^H$. 
Then the pair $(\A^H,m^H)$ is called the \textbf{Ziegler restriction} of $\A$ onto $H$, and the 
map $\pi=\pi_H:D_H(\A) \rightarrow D(\A^H,m^H)$ obtained by taking modulo $\alpha_H$ is called the \textbf{Ziegler restriction map}. The most important results related to 
multiarrangements are the following two.

\begin{thm}[\cite{Z}]
Let $\A$ be free with exponents $(1,d_2,\ldots,d_\ell)$, and 
$H \in \A$. Then $(\A^H,m^H)$ is free with 
exponents $(d_2,\ldots,d_\ell)$ for any $H \in \A$.
\label{Zfree}
\end{thm}

\begin{thm}[\cite{Y2}]
Let $\ell=3$, $H \in \A$ and $\exp(\A^H,m^H)=(d_2,d_3)$. Let 
$\pi:D_H(\A) \rightarrow D(\A^H,m^H)$ be the Ziegler restriction map. Then
$$
\dim_\K \coker \pi=b_2^0(\A)-d_2d_3,
$$
and the equality holds if and only if $\A$ is free with exponents $(1,d_2,d_3)$.
\label{YC}
\end{thm}

If we can determine the whole algebraic structure of $D(\A)$, then of course we can 
see $r(\A)$, which is in general very difficult unless $\A$ is free. By \cite{A5}, we can do it for the arrangement that can be obtained by deleting one hyperplane from free one.

\begin{thm}[Theorem 1.4, \cite{A5}]
Let $\A$ be free with $\exp(\A)=(1,d_2,\ldots,d_\ell)$. Let 
$H \in \A$ with $|\A|-|\A^H|=:d+1$. If $\A':=
\A \setminus \{H\}$ is not free, then $D(\A')$ has a minimal 
set of generators 
$$
\theta_E,\theta_2,\ldots,\theta_\ell,\varphi
$$
such that 
$\theta_E,\theta_2,\ldots,\theta_\ell$ form a basis for 
$D(\A)$ with 
$\deg \theta_i=d_i$, $\varphi \not \in D(\A)$ is of degree $d$, and there is the unique relation 
$$
\sum_{i=1}^\ell f_i \theta_i+\alpha_H\varphi=0
$$
for $\theta_1:=\theta_E,\ f_i \in S$.
\label{free-1}
\end{thm}

\begin{thm}[Theorem 5.5, \cite{A5}]
Let $\ell=3$, $H \in \A$ and $\A':=\A 
\setminus \{H\}$.
Assume that $\A$ is not free. Then $\A'$ is free with $\exp(\A')=(1,d_2,d_3)$ 
if and only if $D(\A)$ is generated by derivations 
$\theta_E,\theta_2,\theta_3,\varphi$ of degrees $\deg \theta_i=d_i+1,\ 
\deg \varphi=|\A^H|-1$ and there is the unique relation 
$$
f_1 \theta_E+f_2\theta_2+f_3 \theta_3+\alpha_H \varphi=0.
$$
\label{free+1}
\end{thm}

Since such arrangements are useful, we give them a name as follows.

\begin{definition}[Definition 1.1, \cite{A5}]
We say that $\A$ is \textbf{plus-one generated (POG)} with 
$\POexp(\A)=(d_1,\ldots,d_\ell)$ and level $d$ if 
$D(\A)$ has a minimal 
set of generators 
$$
\theta_E,\theta_2,\ldots,\theta_\ell,\varphi
$$
such that 
$\deg \theta_i=d_i$, $\deg \varphi=d$, and there is the unique relation 
$$
\sum_{i=1}^\ell f_i \theta_i+\alpha \varphi=0
$$
for some $\alpha \in V^*$. 

We say that $\A$ is \textbf{strongly plus-one generated (SPOG)} with 
$\POexp(\A)=(d_1,\ldots,d_\ell)$ and level $d$ if 
$D(\A)$ has a minimal 
set of generators 
$$
\theta_E,\theta_2,\ldots,\theta_\ell,\varphi
$$
such that 
$\deg \theta_i=d_i$, $\deg \varphi=d$, and there is the unique relation 
$$
\sum_{i=1}^\ell f_i \theta_i+\alpha \varphi=0
$$
for some $\alpha \in V^* \setminus \{0\}$. Such a $\varphi$ is called the \textbf{level element},
and such a set of minimal generators is said to be a \textbf{SPOG-generator}. 
\label{SPOG}
\end{definition}

\begin{rk}
Note that when $\ell=3$ all POG arrangements are 
SPOG by Proposition 5.1, \cite{A5}.
\label{3pogspog}
\end{rk}

\subsection{Addition-deletion theorems on $r(\A)$ for hyperplane arrangements}

First let us show the most fundamental results on 
$r(\A)$. 

\begin{prop}
Let $\ell \ge 2$, $H \in \A,\ 
\A':=\A \setminus \{H\}$. 
Then 
\begin{itemize}
\item[(1)] 
$r(\A') \le r(\A) \le r(\A')+1$. 
\item[(2)]
If $|\A|-|\A^H|>r(\A)$, then $r(\A)=r(\A')$.
\item[(3)]
If $|\A'|-|\A^H|>r(\A')$, then $r(\A)=r(\A')$.
\end{itemize}
\label{fund}
\end{prop}

\proof
Let $L \in \A'$. Since $\mbox{ch}(\K)=0$, it holds that 
$$
D_0(\A') \simeq D_L(\A'),\ 
D_0(\A) \simeq D_L(\A).
$$
So we may consider only derivations in $D_L(\A)$ and 
$D_L(\A')$. 
First assume that $r(\A') \le r(\A)-2$. Let 
$\theta \in D_L(\A')$ be of degree $r(\A')$. Then 
$\alpha_H \theta \in D_L(\A)$ is zero by the definition of $r(\A)$. Thus 
$r(\A')+1 \ge r(\A)$. Since $D_L(\A') \supset D_L(\A)$, $r(\A) \ge r(\A')$, completing the proof of 
(1). 

Next let us prove (2). Assume that $r(\A')+1=r(\A)=:r$. Let 
$0 
\neq \theta \in D_L(\A')_{r-1}$. By Theorem \ref{B}, 
it holds that $D_L(\A')_{<|\A'|-|\A^H|} =D_L(\A)_{<|\A'|-|\A^H|}$. 
In particular, $D_L(\A')_{r-1}=D_L(\A)_{r-1}$. So 
$\theta \in D_L(\A)_{r-1}=(0)$. This is absurd. The same argument shows (3). 
\endproof
\medskip

For an arrangement $\A$, to study $r(\A)$, the Euler derivation does not appear, but it is very important in the following sense.

\begin{lem}
Let $\A \neq \emptyset$ and $ 0 \neq \theta \in D(\A)_d$. Then 
$r(\A) \le d$ if $\theta \not \in S \theta_E$.
\label{fund2}
\end{lem}

\noindent
\proof 
Let 
$\theta':=\theta-\frac{\theta(\alpha_H)}{\alpha_H} \theta_E \in D_H(\A)_d$
for some $H \in \A$. By the assumption, $\theta' \neq 0$, which completes the proof. 
\endproof
\medskip

Now let us introduce the addition-deletion theorems for 
$r(\A)$ in an arbitrary dimension.

\begin{thm}[Addition-deletion theorem for $r(\A)$]
Let $\ell \ge 2$, $H \in \A$ and 
$\A':=\A \setminus \{H\}$. Let $r=r(\A)$, 
$r'=r(\A')$ and
$r''=r(\A''):=r(\A^H)$. 
Then 
$r=r'+1$ if $r'<r''$.
\label{ADgeneral}
\end{thm}

\proof Assume that $r'<r''$ and $r=r'$. Then there is 
$ 0\neq \theta \in D_H(\A')_{r=r'}$ such that $\theta \in D_H(\A)_r$ by Lemma \ref{lemma1}. 
Let $\rho:D(\A) \rightarrow D(\A^H)$ be the Euler 
restriction. Since $r''>r$ and $\deg\rho(\theta)=r<r''$, 
$\rho(\theta)$ is of the form $\rho(f \theta_E)$ by Lemma \ref{fund2}. Since 
$\theta \not \in S \theta_E$, we may replace 
$\theta$ by $\theta - f\theta_E \not \in S \theta_E$ and we may assume that $\rho(\theta)=0$. 
Thus $\theta=\alpha_H \theta'$ with 
$S \theta_E \not \ni\theta' \in D(\A')_{r-1=r'-1}$. By Lemma \ref{fund2}, 
$\theta' \neq 0$ implies that $r(\A') \le r'-1$,  a contradiction.
\endproof

The addition-deletion theorem is related with the restriction theorem 
in general.
For the effect of restriction on $r(\A)$, however, we cannot say  much.

\begin{prop}
Let $\ell \ge 2$, $H \in \A$ and 
$\A':=\A \setminus \{H\}$. Let $r=r(\A)$, 
$r'=r(\A')$ and
$r''=r(\A'')=r(\A^H)$. 
Then
$r'' \le r$ if $r=r'$.
\label{Rgeneral}
\end{prop}

\proof 
By the same proof as in Theorem \ref{ADgeneral}, there is a common 
$0 \neq \theta \in D(\A')_{r=r'} \cap D(\A)_{r}$, and 
we may assume that $0 \neq \rho(\theta)\in D(\A^H)_r \setminus (S/\alpha_HS)\rho(\theta_E) $. Therefore one has $r''\le r$.\endproof

\subsection{The case of free hyperplane arrangements}

We can explicitly describe the behaviour of $r(\A)$ when $\A$ is free.

\begin{thm}
Let $\A$ be free with $\exp(\A)=(1,d_2,\ldots,d_\ell)_\le$. 
Let $H \in \A$ and $\A':=\A \setminus \{H\}$. 
Then $r(\A')=d_2-1$ 
if and only if $|\A|-|\A^H|=d_2$. Otherwise $r(\A')=d_2$. 
\label{freemdr}
\end{thm}

\proof
Recall that $r(\A)=d_2$. 
If $|\A|-|\A^H|=d_2$, then $\A'$ is free with exponents 
$(1,d_2-1,d_3,\ldots,d_\ell)$ by Theorem \ref{MDT}. 
Assume that $r':=r(\A')=d_2-1$. Then $|\A|-|\A^H|$ has to be 
$d_2$ if $\A'$ is free. Assume that $\A'$ 
is not free. Then by Theorem \ref{free-1}, $\A'$ is strictly plus-one generated with exponents 
$(1,d_2,\ldots,d_\ell)_\le$ and level $|\A|-|\A^H|-1=:d$. Let us show that 
$d \ge d_2-1$. Assume that $d < d_2-1$. Note that 
$$
\frac{\overline{Q(\A)}}{\prod_{X \in \A^H} \alpha_{\nu{X} }}\pi(\theta_E)=:\theta_E^H 
\in D(\A^H,m^H)
$$
and 
$d_2>d+1
=\deg \theta_E^H
$,
where $\pi:D_H(\A) \rightarrow D(\A^H,m^H)$ is the Ziegler 
restriction map and $\nu$ is the section in Theorem \ref{B}. 
Since $D(\A^H,m^H)$ is free with exponents $(d_2,\ldots,d_\ell)$ by Theorem \ref{Zfree}, this is a contradiction. 
\endproof

\section{Line arrangements}

\subsection{Addition-deletion theorems on $r(\A)$ for line arrangements}

First, let us recall Terao's 
addition-deletion theorem for line arrangements. 

\begin{thm}[Terao's addition-deletion theorem, \cite{T1}]
Let $\ell=3$, $H \in \A$ and $\A':=\A 
\setminus \{H\}$. Then 

(1)\,\,
$\A$ is free with $\exp(\A)=(1,a,b+1)$ if $\A'$ is free with $\exp(\A')=(1,a,b)$ and 
$|\A^H|=a+1$. 

(2)\,\,
$\A'$ is free with $\exp(\A')=(1,a,b-1)$ if $\A$ is free with $\exp(\A)=(1,a,b)$ and 
$|\A^H|=a+1$. 

\label{adddel3}
\end{thm}

When $\ell=3$, Theorem \ref{ADgeneral} is more combinatorial. 

\begin{thm}[Addition theorem for $r(\A)$]
Let $\ell=3$ and $\A:=\A' \cup \{H\}$ with $H \not \in \A'$. Assume that 
$r(\A')=r'$. If 
$|\A^H| \ge r'+2$, then $r(\A)=r'+1$. 
\label{main}
\end{thm}

\proof 
Since $\exp(\A^H)=(1,|\A^H|-1)$, $|\A^H| \ge r'+2$ shows that 
$r''>r'$. Now apply Theorem 
\ref{ADgeneral}. \endproof

\begin{thm}[Deletion theorem for $r(\A)$]
Let $\ell=3$ and $\A:=\A' \cup \{H\}$ with $H \not \in \A'$. Assume that 
$r(\A)=r$. If 
$|\A^H| \ge r+2$, 
then $r(\A')=r-1$. 
\label{main2}
\end{thm}

\proof
Apply the same proof as in Theorem \ref{main}. \endproof

\subsection{The case of free line arrangements}

By the same reason as above, we can describe $r(\A)$ when $\ell=3$ and 
the arrangement $\A$ is free. 

\begin{thm}
Let $\ell=3$, $H \in \A$ and $\A':=\A \setminus \{H\}$. 

(1)\,\,
Assume that $\A$ is free with $\exp(\A)=(1,d_2,d_3)_\le$. Then $r(\A')=d_2-1$ 
if and only if $|\A^H|=1+d_3$. Otherwise $r(\A')=d_2$.

(2)\,\,
Assume that $\A'$ is free with $\exp(\A')=(1,d_2,d_3)_\le$. Then $r(\A)=d_2+1$ 
if and only if $|\A^H|\neq 1+d_2$. Otherwise $r(\A')=d_2$.
\label{freemdr3}
\end{thm}

\proof (1) follows from 
Theorem \ref{freemdr}. Let us show (2). 
By \cite{A}, $\A$ is free with exponents $(1,d_2+1,d_3)$ if and only if $
|\A'|-|\A^H|=d_3$, and 
$$
|\A'|-|\A^H|=d_3,\ \mbox{or}\ 
|\A'|-|\A^H|\le d_2.
$$
If $\A$ is not free, then $\A$ is strictly plus-one generated with exponents 
$(1,d_2+1,d_3+1)$ and level $d=|\A^H|-1$ by Theorem \ref{free+1}. Since $\ell =3$, $\A$ is SPOG. Thus 
$d \ge d_3+1$. Thus $mdr(A)=d_2+1$. Since 
this occurs if and only if $|\A'|-|\A^H| \neq d_2,d_3$, it suffices to show that 
$r:=mdr(\A)=d_2$ if $|\A'|-|\A^H| = d_3$, which is trivial since $\exp(\A)=(1,d_2,d_3+1)$ by Theorem \ref{adddel3}. \endproof

\subsection{Points of high multiplicity and the invariant $r(\A)$}
In this subsection $\A':f'=0$ is a line arrangement, and $p=(1:0:0)$ is an intersection  point
on $\A'$ of maximal multiplicity, say $m'=\mult(\A',p)$. To this situation, one can associate a primitive Jacobian syzygy as explained in
\cite[Section 2.2]{DMich}. We recall this construction here.
Let $g=0$ be the equation of the subarrangement of $\A'$ formed by the $m'$ lines in $\A'$ passing through $p$. Then we can write $f'=gh$ for some polynomial $h \in S$. Since $g$ is a product of linear factors of the form $sy+tz$, it follows that $f'_x=gh_x$ and hence $g=G.C.D.(f',f'_x)$. The syzygy constructed as explained there is primitive
and has degree $r'_p=d'-m'$.
As shown in \cite[Theorem 1.2]{DMich}, the following cases are possible for $r'=r(\A')$.

\medskip

\noindent {\bf Case A:} $r'=r'_p=d'-m'$, in other words the constructed syzygy has minimal degree.

\medskip 
\noindent {\bf Case B:} $r'<r'_p=d'-m'$, in other words the constructed syzygy has not minimal degree. Then two situations are possible, namely
\medskip 

\noindent {\bf Subcase B1:} $r'=m'-1$, and then $2m'<d'+1$ and $\A'$ is free
with exponents $d_1=1, d_2=m'-1<d_3=d'-m'$, or 

\medskip 

\noindent {\bf Subcase B2:} $m' \leq r' \leq d'-m'-1$, and then $2m' <d'$.

\medskip 

\noindent This discussion implies the following.
\begin{cor}
\label{cor1}
If the line arrangement $\A'$ satisfies $2m' \geq d'$, then $r'=r(\A')$ is determined by the weak combinatorics of $\A'$.
\end{cor}

\proof
If $2m' \geq d'+1$, it follows that only Case A is possible, and hence $r'=d'-m'$. When $2m'=d'$, then both Case A and Subcase B1 are possible, hence we have either $r'=m'-1$ or $r'=d'-m'=m'$.
If $\A'$ is in the situation of Subcase B1, then we know that
$$\tau(\A')=(d'-1)^2-r'(d-r'-1)=(d'-1)^2-(m'-1)(d'-m').$$
On the other hand, if $\A'$ is in the situation of Case A, then we know that
$$\tau(\A') \leq (d'-1)^2-r'(d-r'-1)-1=(d'-1)^2-(d'-m')(m'-1)-1,$$
see \cite{dPW,Dmax}. Since the total Tjurina number is determined by the weak combinatorics, recall \eqref{TauW},
 this completes the proof.
\endproof
\begin{rk}
\label{rk1.5}
In Ziegler's celebrated example, see \cite{Zi}, we have two line arrangements $\A'_1$ and $\A'_2$ of degree $d'=9$ and such that $m'=3$ in both cases.
For one of them, say for $\A'_1$,  the six triple points are on a conic, and one has $r'_1=5 =d'-m'-1$, hence we are in Subcase B2 above.
For the other one, say for $\A'_2$, the six triple points are not on a conic,
and one has $r'_2=6=d'-m'$, so we are in Case A. This shows that the combinatorics of $\A'$ cannot decide in which case {\bf  A}, {\bf  B1} or {\bf  B2} we are in the above discussion. This example is discussed in \cite[Remark 8.5]{DHA}.
One can find there some equations for the arrangements $\A'_1$ and $\A'_2$, namely
$$\A'_1: xy(x-y-z)(x-y+z)(2x+y-2z)\times$$
$$\times (x+3y-3z)(3x+2y+3z)(x+5y+5z)(7x-4y-z)=0$$
 and
$$\A'_2: xy(4x-5y-5z)(x-y+z)(16x+13y-20z)\times$$
$$\times(x+3y-3z)(3x+2y+3z)(x+5y+5z)(7x-4y-z)=0.$$
In fact, the equation for $\A'_2$ given in \cite[Remark 8.5]{DHA} is not correct, and we take the opportunity here to correct this equation.
\end{rk}

\section{Applications}

\subsection{Line arrangements which are Tjurina maximal}
Recall that the {\it global Tjurina number} $\tau(C)$ of the plane curve $C:f=0$ can be defined as either the degree of the Jacobian ideal $J_f=(f_x,f_y,f_z)$, or as the sum of the Tjurina numbers of all the singularities of the curve $C$. It was shown by A. du Plessis and C.T.C. Wall that one has the following result, see \cite[Theorem 3.2]{dPW}, and also \cite[Theorem 20]{E} for a new approach.

\begin{thm}
\label{thmCTC}
Let $C:f=0$ be a reduced plane curve of degree $d$ and let $r=mdr(C)$.
Then the following hold.
\begin{enumerate}

\item If $r <d/2$, then $\tau(C) \leq \tau(d,r)_{max}= (d-1)(d-r-1)+r^2$ and the equality holds if and only if the curve $C$ is free.

\item If $d/2 \leq r \leq d-1$, then
$\tau(C) \leq \tau(d,r)_{max}$,
where, in this case, we set
$$\tau(d,r)_{max}=(d-1)(d-r-1)+r^2-{ 2r-d+2 \choose 2}.$$

\end{enumerate}

\end{thm}
The curve $C:f=0$ in this Theorem is called {\bf maximal Tjurina of type $(d,r)$} if one has the equality
$$\tau(C) = \tau(d,r)_{max}.$$
The characterization and the existence of maximal Tjurina curves of type $(d,r)$, with $d/2 \leq r \leq d-1$ is discussed in \cite{DS}. 
In this note we prove the existence of maximal Tjurina curves of type $(d,r)$ in many cases.  We start with the following.

\begin{prop}
\label{prop100}
If $\A':f'=0$ is a Tjurina maximal line arrangement of type $(d',r')$ with $r' \geq (d'-1)/2$ and $H$ is a new line in $\PP^2$ such that
$$|\A' \cap H|=r'+2,$$
then $\A=\A' \cup H$ is a Tjurina maximal line arrangement of type $(d,r)$ with $d=d'+1$ and $r=r'+1$.
\end{prop}
\proof
First note that Theorem \ref{main} implies that $r=r(\A)=r(\A')+1=r'+1$. Hence to show that $\A$ is a Tjurina maximal line arrangement of type $(d,r)$, it is enough to show that it has the global Tjurina number
$\tau(\A)=\tau(d,r)_{max}$, given by the formula in Theorem \ref{thmCTC} (2). A direct computation shows that this is equivalent to  the following
\begin{equation}
\label{e100}
\tau(\A)-\tau(\A')=2d'-r'-2.
\end{equation}
To measure the difference $\tau(\A)-\tau(\A')$, assume that $\A' \cap H$ consists of $s$ points, say $p_1,\ldots,p_s$, with multiplicities
$m_1, \ldots,m_s$ regarded as points on $\A'$. 
When we add the new line $H$, the point $p_j$ will have multiplicity $m_j+1$, so the increase in Tjurina number at $p_j$ is 
$$m_j^2-(m_j-1)^2=2m_j-1.$$
It follows that
\begin{equation}
\label{e200}
\tau(\A)-\tau(\A')=\sum_{j=1,s}(2m_j-1)=2d'-s.
\end{equation}
This ends the proof of the claim.
\endproof

\begin{thm}
\label{thm100}
Given a pair of positive integers $(d,r)$ such that $d\geq 4$ and
$$\frac{d}{2} \leq r \leq \frac{2}{3}(d-1),$$
then there is a real line arrangement $\A$ in $\PP^2$ which is Tjurina maximal  of type $(d,r)$.
\end{thm}
\proof
We set $r=d-k$ for some $k \geq 2$, and the equalities involving $d$ and $r$ in Theorem
\ref{thm100} are equivalent to
$$2k \leq d \leq 3k-2.$$
Hence we have to show the existence of a real line arrangement $\A$ in $\PP^2$ which is Tjurina maximal  of type $(d,d-k)$, where $d$ and $k$ satisfy the above inequalities.
We start with the line arrangement
$$\A_0: f_0(x,y,z)=x(x-z) \ldots (x-(k-2)z)y(y-z) \ldots (y-(k-2)z)z=0,$$
which is free, even supersolvable, and also 
Tjurina maximal  of type $(2k-1,k-1)$.
If we add the line
$$H_1: y=x+z,$$
and apply Proposition \ref{prop100} with $\A'=\A_0$ and $H=H_1$. We get that
$$\A_1=\A_0\cup H_1: f_1(x,y,z)=f_0(x,y,z)(x-y+z)=0$$
is a Tjurina maximal line arrangement of type $(2k,k)$.
Then we add the line
$$H_2: y=x+2z,$$
and apply Proposition \ref{prop100} with $\A'=\A_1$ and $H=H_2$. We get that
$$\A_2=\A_1\cup H_2: f_2(x,y,z)=f_0(x,y,z)(x-y+z)(x-y+2z)=0$$
is a Tjurina maximal line arrangement of type $(2k+1,k+1)$.
Assume now that $\A_j$ has been constructed,
for $2 \leq j < k-2$, and it is a Tjurina maximal line arrangement of type $(2k+j-1,k+j-1)$. Then we construct $\A_{j+1}$ by adding the new line
$$H_{j+1}: y=x+(j+1)z$$
and apply Proposition \ref{prop100} with $\A'=\A_j$ and $H=H_{j+1}$. We get that
$$\A_{j+1}=\A_j\cup H_{j+1}: f_{j+1}(x,y,z)=f_0(x,y,z)(x-y+z)\ldots (x-y+(j+1)z)=0$$
is a Tjurina maximal line arrangement of type $(2k+j,k+j)$.
This construction ends when we construct $\A_{k-1}$, because after this value the hypothesis of Proposition \ref{prop100} is no longer verified.
\endproof

\begin{cor} 
\label{corex5.5}
For any odd degree $d=2r-1\geq 7$, there is a maximal Tjurina real line arrangement of type $(2r-1,r)$.
\end{cor} 
\proof
Just consider the arrangement $\A_2$ in the above proof.
\endproof
Note that the last arrangement $\A_{k-1}$ constructed in the proof of Theorem \ref{thm100} consists of the line at infinity $z=0$ and three families of parallel lines, each containing $k-1$ lines. Hence this arrangement has 3 points of maximal multiplicity equal to $k$ on the line at infinity. 
When $d-r=k=2k'+1$ is odd, we can continue the above construction and get a stronger result.
\begin{prop}
\label{prop101}
Given a pair of positive integer $(d,r)$ such that $d\geq 4$, $k=d-r$ is odd and
$$\frac{d}{2} \leq r \leq \frac{3}{4}(d-1),$$
then there is a real line arrangement $\A$ in $\PP^2$ which is Tjurina maximal  of type $(d,r)$.
\end{prop}

\proof When $d-r=k=2k'+1$ is odd, we can continue the above construction in two steps, as follows.
To get $\A_{k}$ from $\A_{k-1}$ we add the line
$$H_m: x+y=3k'z.$$
Using Proposition \ref{prop100} we get that $\A_k$ is  a Tjurina maximal line arrangement of type $(3k-1,2k-1)$. Then, in the first step, we  add the lines
$$H_j: x+y=(3k'+j-k)z,$$
for $j=2k'+2,..., 3k'$, to get new line arrangements
$\A_j=\A_{k-1} \cup H_k \cup \ldots \cup H_j$,
which are  Tjurina maximal  of type $(2k+j-1, k+j-1)$ for each $j=2k'+2,..., 3k'$. If we increase the coefficient of $z$ beyond this value
$6k'-k=2k-3$, the number of intersection points in
$\A' \cap H$ is no longer a strictly increasing sequence $q, q+1,q+2,...$ as until now, but has repetitions of the form $q,q, q+1,q+1,q+2,q+2,...$, and hence we have a choice in selecting the new line to add between two possibilities.
This is the second step in this construction.
The largest type we can get in this way is $(4k-3,3k-3)$, and we denote such an arrangement by $\A_{2k-2}$, since it is obtained from $\A_0$ by adding $2k-2$ lines. 
\endproof

Note that the last arrangement $\A_{2k-2}$ constructed in the proof of Proposition \ref{prop101} consists of the line at infinity $z=0$ and four families of parallel lines, each containing $k-1$ lines. Hence this arrangement has 4 points of maximal multiplicity equal to $k$ on the line at infinity. 
\begin{ex}
\label{ex100}
As an illustration, consider the case $k=7$, and hence $k'=3$. The first sequence of line arrangements constructed in the proof of Theorem \ref{thm100} has the following equations
$$\A_0: f_0(x,y,z)=x(x-z) \ldots (x-5z)y(y-z) \ldots (y-5z)z=0,$$
and
$$\A_j: f_j(x,y,z)=f_0(x,y,z)(x-y+z)(x-y+2z) \ldots (x-y+jz)=0,$$
for $j=1,2,...,6$.
The  arrangements constructed in the first step in Proposition \ref{prop101}
consists of the following two arrangements
$$\A_7: f_7(x,y,z)=f_6(x,y,z)(x+y-9z)\text{ and }  \A_8: f_8(x,y,z)=f_7(x,y,z)(x+y-10z).$$
The second sequence of arrangements constructed in Proposition \ref{prop101}, when a choice is possible,
consists of the following four arrangements obtained by taking the minimal absolute value for the coefficient of $z$:
$$\A_9: f_{9}(x,y,z)=f_8(x,y,z)(x+y-11z),  \  \
\A_{10}: f_{10}(x,y,z)=f_9(x,y,z)(x+y-13z),$$
$$\A_{11}: f_{11}(x,y,z)=f_{10}(x,y,z)(x+y-15z),$$
and
$$\A_{12}: f_{12}(x,y,z)=f_{11}(x,y,z)(x+y-17z).$$
Note that $\A_0$ has $2k-1=13$ lines, and $\A_{12}$ has $13+12=25=4k-3$ lines, as expected.
Here $\A_j$ is a Tjurina maximal line arrangement of type $(13+j,6+j)$
for $j=0,1,...,12$.
\end{ex}

\begin{rk}
\label{rk101}
When $k=2k'$ is even, then there are two cases.
When $k'$ is odd, we have found no simple way to add a new line to $\A_{k-1}$ in order to get a larger Tjurina maximal line arrangement.
This is due to the fact that the number of intersection points in
$\A' \cap H$ in this case, for $H$ a line of the form to $x+y-az=0$,  has repetitions and gaps of the form $2q+1,2q+1,2q+3,2q+3,...$, i.e. we get only odd numbers, and hence it is impossible to apply Proposition \ref{prop100}.
When $k'$ is even, we can construct $\A_k$ as follows: we add the line
$$H_k: x+y=(3k'-2)z.$$
Using Proposition \ref{prop100} we get that $\A_k$ is  a Tjurina maximal line arrangement of type $(3k-1,2k-1)$.
However, this construction stops here, since in this case the number of intersection points in
$\A' \cap H$ has repetitions and gaps of the form $2q,2q,2q+2,2q+2,...$.
\end{rk}
\begin{ex}
\label{ex101}
As an illustration, consider the case $k=8$, and hence $k'=4$. The first sequence of line arrangements constructed in the proof of Theorem \ref{thm100} has the following equations
$$\A_0: f_0(x,y,z)=x(x-z) \ldots (x-6z)y(y-z) \ldots (y-6z)z=0,$$
and
$$\A_j: f_j(x,y,z)=f_0(x,y,z)(x-y+z)(x-y+2z) \ldots (x-y+jz)=0,$$
for $j=1,2,...,7$.
The largest line arrangement constructed in Remark \ref{rk101} is the following:
$$\A_8: f_8(x,y,z)=f_7(x,y,z)(x+y-10z).$$
Here $\A_j$ is a Tjurina maximal line arrangement of type $(15+j,7+j)$
for $j=0,1,...,8$.
\end{ex}

Consider now the question of the existence of maximal Tjurina line arrangements $\A$ of $d$ lines with large invariant $r=r(\A)$ with respect to $d$. 
The first case we consider is  $r=d-4$, see Remark \ref{rkSE} below.
To give a positive answer to this question, consider the following arrangements 
$\A_{3p+2}$ for $p\geq 2$, defined by 
$$
xy(\frac{x}{2^{p+1}}+\frac{y}{3^{p+1}}-z)\prod_{j=1}^{p} (\frac{x}{2^j}+\frac{y}{3^j}-z)
(\frac{x}{2^j}+\frac{y}{3^{j+1}}-z) 
\prod_{j=1}^{p-1} (\frac{x}{2^j}+\frac{y}{3^{j+2}}-z)=0,
$$
$$
\A_{3p+3}=\A_{3p+2} \cup \{H_1\}
$$
with $H_1:27x-8y=0$, and 
$$
\A_{3p+4}=\A_{3p+3} \cup \{H_2\}
$$
with $H_2:x-y=0$. Using these three families of line arrangements, we can define $\A_d$ for all $d\geq 8$. 
The following result proves a conjecture made in \cite{DS}, where it was shown that $\A_d$ is a maximal Tjurina line arrangement of type $(d,d-4)$ if and only if $r(\A_d)=d-4$.

\begin{thm}
Let $\A_d$ be the arrangement defined above. Then $r(\A_{d})=d-4$ for all $d \geq 8$. 
\label{d-4}
\end{thm}

\proof
We know that the statement is true when $d$ is small by \cite{DS}. We apply Theorem \ref{main} repeatedly in the following addition steps:
\begin{eqnarray*}
\A_{3p+2} &\rightarrow&  \A_{3p+5} \\
\A_{3p+2} &\rightarrow&  \A_{3p+3}\\
\A_{3p+3} &\rightarrow&  \A_{3p+4}.
\end{eqnarray*}
In fact, in the first step, we add 3 lines, so we apply Theorem \ref{main}
three times.
Then an elementary counting of the intersection points completes the proof. \endproof

Finally, consider  the question of the existence of maximal Tjurina line arrangement of $d$ lines with  $r=d-3$.
Define the families of line arrangements
$\B_{2p}$ by 
$$
xy\prod_{j=1}^{p-1} (\frac{x}{2^j}+\frac{y}{3^j}-z)
(\frac{x}{2^j}+\frac{y}{3^{j+1}}-z) =0,
$$
and $\B_{2p+1}$ by 
$$
xy (\frac{x}{2^d}+\frac{y}{3^d}-z)\prod_{j=1}^{d-1} (\frac{x}{2^j}+\frac{y}{3^j}-z)
(\frac{x}{2^j}+\frac{y}{3^{j+1}}-z) =0,
$$
Using these two families, we can define $\B_d$ for all $d\geq 7$. In \cite{DS} it was shown that $\B_d$ is a maximal Tjurina line arrangement of type $(d,d-3)$ if and only if $r(\B_d)=d-3$, and it was conjectured that one has $r(\B_d)=d-3$.
 By the same proof as in Theorem \ref{d-4}, we can show it.

\begin{thm}
Let $\B_d$ be the arrangement defined above. Then $r(\B_{d})=d-3$ for all $d \geq 7$. 
\label{d-3}
\end{thm}

\begin{rk}
\label{rkSE}
For any line arrangement with $d \geq 2$, one has $r \leq d-2$, by our discussion before Corollary \ref{cor1}.
 Any generic line arrangement is maximal Tjurina of type $(d,d-2)$ when $d \geq 2$,  see \cite[Remark 2.2]{DS}, hence the cases $r=d-4$ and $r=d-3$ considered in Theorems \ref{d-4} and \ref{d-3} are the largest possible values of $r$ where the existence of a 
maximal Tjurina line arrangement of type $(d,r)$ is not obvious.
\end{rk}

\subsection{Adding a generic line to a line arrangement}
First we take a generic secant passing through the point $p$ of maximal multiplicity $m'$ for the line arrangement $\A'$.
\begin{prop}
\label{prop2}
Let $\A':f'=0$ be a line arrangement, let $L$ be a generic line passing through the maximal multiplicity intersection point $p \in \A'$, and let $\A=\A' \cup L:f=0$.
Then one has the following.
\begin{enumerate}
\item The (weak) combinatorics of the line arrangement $\A'$ determines the (weak) combinatorics of the line arrangement $\A$. In particular, one has
$$\tau(\A)=\tau(\A')+d'+m'-1.$$
\item $r= r'+1$ if $r' <d'-m'$. 
\item  $r=r'$ if $r'=d'-m'$.
\end{enumerate}

\end{prop}
\proof
The first claim is obvious. To prove the second claim, we use Theorem \ref{main}. The number of intersection points of $\A$ on $L$
is $|\A^L|=1+d'-m'$ and hence the  condition in Theorem \ref{main}, namely 
$$|\A^L| \geq r'+2$$
is equivalent to
$$1+d'-m'\geq r'+2$$
or
$$r' \leq d'-m'-1.$$
This condition is satisfied by our assumption $r'<d'-m'$.
To prove the claim (3), note that by Proposition \ref{fund} we have
$d'-m'=r' \leq r$.
On the other hand, \cite[Theorem 1.2]{DMich} implies
$r \leq d-m=(d'+1)-(m'+1)=d'-m'$, where $m=\mult(\A,p)$.
\endproof
Next we add a generic line, meeting $\A'$ only at simple points.

\begin{prop}
\label{prop4}
Let $\A':f'=0$ be a line arrangement with $d' \geq 2$ and  $L$ be a generic line. Consider the new line arrangement $\A=\A' \cup L:f=0$.
Then one has the following.
\begin{enumerate}
\item The (weak) combinatorics of the line arrangement $\A'$ determines the (weak) combinatorics of the line arrangement $\A$. In particular, one has
$$\tau(\A)=\tau(\A')+d'.$$
\item $r= r'+1$.
\end{enumerate}
\end{prop}
\proof
The intersection points in $\A$ which are not the same as the corresponding ones in $\A'$ are the $d'$ double points along the line $L$.
These points add $d'$ to the global Tjurina number of $\A$.
 The number of intersection points of $\A$ on $L$
is $|\A^L|=d'$ and hence the condition in Theorem \ref{main}, namely 
$$|\A^L| \geq r'+2$$
is equivalent to
$$d' \geq r'+2$$
This condition is satisfied, since $r' \leq d'-m' \leq d'-2$.

\endproof
\begin{ex}
\label{ex3}
Consider again  Ziegler's arrangements $\A'_1$ and $\A'_2$ from Remark 
\ref{rk1.5}. If we apply Proposition \ref{prop2} (2) to the arrangement $\A'_1$, we get a new arrangement with $r_1=r'_1+1=5+1=6$.
If we apply Proposition \ref{prop2} (3) to the arrangement $\A'_2$, we get a new arrangement with $r_2=r'_2=6$. Hence by adding a generic line through a triple point, the difference between 
$r_1'$ and $r_2'$ disappears.
On the other hand, if we add a generic line $L$ to both arrangements $\A'_1$ and $\A'_2$, we get again two line arrangements $\A_1$ and $\A_2$ with $d=10$, $r_1=5+1=6 <r_2=6+1=7$ and having the same combinatorics. By continuing to add generic lines we can construct such pairs for any $d \geq 9$.
\end{ex}

\subsection{On line arrangements with double and triple points}
Note that Ziegler's arrangements $\A'_1$ and $\A'_2$ have both only double and triple points, more precisely $n_2= 18$ and $n_3=6$.
The following result says that $n_3=6$ is the minimal value for which such pairs with the same combinatorics but distinct values for $r$ can be constructed.
\begin{thm}
\label{thm2}
Let $\A$ be a line arrangement with $d=|\A| \geq 2$, having $n_2$ double points, $n_3$ triple points and no points of higher multiplicity.
\begin{enumerate}
\item If $n_3=0$, then $r(\A)=d-2$.

\item If $1 \leq n_3 \leq 3$, then $r(\A)=d-3$.

\item If $n_3=4$, then $r(\A)=d-3$, unless any line of the arrangement $\A$, passing through a triple point of $\A$, contains an extra triple point of $\A$. In this latter situation,
$\A$ is obtained, up-to a change of coordinates, from the
arrangement 
$$\A(2,2,3): (x^2-y^2)(x^2-z^2)(y^2-z^2)=0$$
by adding $d-6$ generic lines, and then $r(\A)=d-4$.

\item If $n_3=5$, then there are two possibilities.

\subitem (A) There is at least one triple point $p$ in $\A$
and a line $L$ in $\A$,
passing through $p$ and containing only $p$ as a triple point. If the line arrangement $\A'=\A \setminus L$ is obtained, up-to a change of coordinates, from the arrangement
$\A(2,2,3)$ by adding $d-7$ generic lines, then
$r(\A)=d-4$. Otherwise $r(\A)=d-3$.

\subitem (B) For any of the 5 triple points, the 3 lines meeting at this point contain each at least an extra triple point, and then $r(\A)=d-4$. The intersection lattice of $\A$ in this case is the same as the intersection lattice of the arrangement obtained by adding $d-7$ generic lines to the following arrangement
$$\B: y(y+x)(y-x)(y+x-2z)(y-x-2z)(3y+x-2z)(3y-x-2z)=0.$$

\end{enumerate}

\end{thm}
\proof
The claim (1) is well known, see for instance \cite[Theorem 4.1] {DStEdin}.

Consider now the claim (2). Let $p$ be a triple point, and note that, since $n_3 \leq 3$, there is a line $L$ in $\A$, passing through $p$, containing only $p$ as a triple point. Since $r(\A)\leq d-m=d-3$, it is enough to show that $r=r(\A) \leq d-4$ leads to a contradiction.

Apply Theorem \ref{main2} to the arrangement $\A'=\A\setminus L$ and the line $L$.
 The number of intersection points of $\A$ on $L$
is $|\A^L|=d-2$ and hence the condition in Theorem \ref{main2}, namely 
$$|\A^L| \geq r'+2,$$
where $r'=r(\A')$, is satisfied by our assumption.
It follows that $r'=r-1 \leq d-5=d'-4$.
We start with the case $n_3=1$. Then $\A'$ is nodal, so $r'=d'-2$ a contradiction. Hence in this case $r=d-3$. The cases $n_3=2$ and $n_3=3$ can be treated in exactly the same way, using the previous
cases.

To treat the claim (3), note that there are two possibilities.
The first one is that there is a triple point $p$ and a line $L$ in $\A$,
passing through $p$ and containing only $p$ as a triple point. Then we can repeat the argument in the case (2) and get $r=d-3$.
The second case is when, for any of the 4 triple points, the 3 lines meeting at this point contain each an extra triple point. This situation occurs for the arrangement $\A(2,2,3)$, and it is known that this arrangement has $r=2=d-4$. If we are in this situation, the 6 lines determined by the 4 triple points form an arrangement which is, up-to a linear change of coordinates, the arrangement $\A(2,2,3)$. The additional lines must create only double points, so they are generic lines.
Using Proposition \ref{prop4}, we see that for any arrangement $\A$ constructed in this way we get $r=d-4$.

In the final claim (4), if we are in case (A), we can delete the line $L$, and the resulting arrangement $\A'$ has $n_3=4$. Hence the two cases discussed in (3) are possible.
More precisely, we know that $r \leq d-m=d-3$. Assume $r \leq d-4$ and apply Theorem \ref{main2}. We get as above
$r' =r-1\leq d-5=d'-4$. Using the claim (3), we infer that in this case 
$\A'$ is obtained from the arrangement
$\A(2,2,3)$ by adding $d-7$ generic lines.

If we are in case (B), it is enough to check that for the arrangement $\B$ of 7 lines we have $r(\B)=3$, which follows by a direct computation using SINGULAR, and then we use Proposition \ref{prop4}.

The possible configurations of the 5 triple points in $\A$ are discussed next, and this discussion shows that only the situations (A) and (B) are possible.

\medskip

\noindent {\bf Case 1:} Assume first that each line in $\A$ contains at most 2 triple points. If each triple point $p$ is connected to 3 other triple points by lines in $\A$, it means that there is a unique triple point $p'$ not connected to $p$. The 5 triple points are in this way divided in a number of pairs $\{p,p'\}$, 
a contradiction. Hence  in this case we are in the situation (A).

\medskip

\noindent {\bf Case 2:} Assume next that there is a unique line $L'$ in $\A$ containing 3 or more triple points. If $L'$ contains at least 4 triple points, the claim is clear, any triple point $p$ on $L'$ is a good choice, to see that we are again in situation (A). Assume now that $L'$ contain 3 points $p_1,p_2$ and $p_3$, and the remaining triple points are $q_1$ and $q_2$ are not on $L'$. Each of the points $p_j$ has to be connected with both points $q_1$ and $q_2$, in order to avoid being again in the situation (A). In this way, by considering these 7 lines,  we get a line arrangement with the same  combinatorics as $\B$. Therefore, if this happens, we are in the situation (B).

\medskip

\noindent {\bf Case 3:} Assume finally that there are two lines $L'$ and $L''$ in $\A$ containing each 3 triple points. Then the intersection point $p=L' \cap L''$ has to be a triple point, and the third line through $p$, call it $L$, contains no triple points except $p$. Therefore we are in the situation (A).
\endproof

\begin{ex}
\label{ex10}
(i) The arrangement 
$$\A: f=xyz(x+z)(x+y-z)(2x-y)(x+y+3z)=0$$
has $d=7$, $n_3=4$ and $r=d-3=4$. Hence it illustrates the case (3), when $r=d-3$ in Theorem \ref{thm2}. The line $L:x+y-z=0$ contains a unique triple point, namely the point $p=(1:-1:0)$.

\medskip

\noindent (ii) The arrangement 
$$\A: f=xyz(x+z)(x+y-z)(2x-y)(x+y+3z)(y+z)=0$$
has $d=8$, $n_3=5$ and $r=d-3=5$. Hence it illustrates the case (4), subcase (A), when $r=d-3$ in Theorem \ref{thm2}. The line $L:x+y-z=0$ contains a unique triple point, namely $p=(1:-1:0)$, and also the line
$L':y+z=0$ contains a unique triple point, namely $p'=(1:0:0)$.

\end{ex}

\begin{cor}
\label{cor2}
Let $\A:f=0$ be a line arrangement with $d=|\A| \geq 2$, having $n_2$ double points, $n_3$ triple points and no points of higher multiplicity.
The invariant $r(\A)$ is determined by the combinatorics of $\A$ if and only if $n_3 \leq 5$.
\end{cor}
\proof
The claim that $r(\A)$ is determined by the combinatorics of $\A$ if $n_3 \leq 5$ follows from Theorem \ref{thm2}.
The claim that $r(\A)$ is not determined by the combinatorics of $\A$ if $n_3 =6$ follows from Ziegler's example of arrangements $\A'_j$  for $j=1,2$ discussed in Remark \ref{rk1.5} and in Example \ref{ex3}, where we show that any $d \geq 9$ can be realized. To increase the number of triple points, it is enough to pick a double point $p$ in $\A'_j$  for $j=1,2$ and add a generic line $L$ passing through $p$. We can apply Theorem  \ref{main} and show that the arrangements $\A_j=\A'_j \cup L$ have
$r_1=r_1'+1=6$ and $r_2=r_2'+1=7$, and they both have $n_3=7$ and the same combinatorics. Note that $|\A_1^L|=8>r_1'+2=7$, and
$|\A_2^L|=8=r_2'+2$, hence we need the full strength of Theorem  \ref{main}.

Proceeding in this way, it is clear that for any $n_3 \geq 6$, one can construct a pair of line arrangements having only double points and $n_3$ triple points, with the same combinatorics, but distinct invariants $r$.
\endproof

\begin{rk}
\label{rkNU}
Using \cite[Proposition 3.2, (3) and (4)]{AD} and \cite[Theorem 3.2 (1)]{DStRLM}, it follows that Conjecture \ref{conjT4} holds for the line arrangements having only double and triple points. More pecisely, \cite[Theorem 3.2 (1)]{DStRLM} shows that an arrangement  $\A$ of $d$ lines, having only double and triple points, satisfy
$$r(\A) \geq \frac{d-2}{2}.$$
Then \cite[Proposition 3.2, (3) and (4)]{AD} shows that in this case, the generic splitting type $(e_1,e_2)$ is determined by $d$.
\end{rk}


\begin{thebibliography}{00}
\bibitem{A}
T. Abe, 
Roots of characteristic polynomials and 
and intersection points of line arrangements. 
\textit{J. Singularities}, 
\textbf{8} (2014), 100--117.

\bibitem{A2}
T. Abe,
Divisionally free arrangements of hyperplanes. 
\textit{Invent. Math.} \textbf{204} (2016), no. 1, 317--346.



\bibitem{A4}
T. Abe, 
Deletion theorem and combinatorics of hyperplane arrangements.
\textit{Math. Ann.} \textbf{373} (2019), issue 1-2, 581--595. 


\bibitem{A5}
T. Abe, 
Plus-one generated and next to free arrangements of hyperplanes, 
\textit{Int. Math. Res. Not.}.
DOI:10.1093/imrn/rnz099


\bibitem{AD}
T. Abe and A. Dimca, 
On the splitting types of bundles of logarithmic vector fields along plane curves. 
\textit{Int. J. Math}.,  \textbf{29}, no. 8, 1850055, 20 pp.


\bibitem{AT}
T. Abe and H. Terao,
Multiple addition, deletion and restriction theorems for hyperplane arrangements. 
\textit{Proc. Amer. Math. Soc.}. DOI:10.1090/proc/14592










\bibitem{CHMN} D. Cook, B. Harbourne, J. Migliore, U. Nagel,  Line arrangements and configurations of points with an unexpected geometric property. \textit{Compositio Math.} \textbf{154}(2018), 2150--2194.



\bibitem{DHA}  A. Dimca,   {\em Hyperplane Arrangements: An Introduction}, Universitext, Springer, 2017

\bibitem{DMich}  A. Dimca, Curve arrangements, pencils, and Jacobian syzygies,  Michigan Math. J. 66 (2017), 347--365.


\bibitem{Dmax}  A. Dimca, Freeness versus maximal global Tjurina number for plane curves, 
Math. Proc. Cambridge Phil. Soc.  163 (2017), 161--172.

\bibitem{Drcc} A. Dimca, On rational cuspidal plane curves, and the local cohomology of Jacobian rings, arXiv:1707.05258.  to appear in Commentarii Mathematici Helvetici.


\bibitem{DIM} A. Dimca, D. Ibadula, A. M\u acinic, Freeness for 13 lines arrangements is combinatorial, Discrete Mathematics 342 (2019), 2445--2453. 








\bibitem{DStEdin} A. Dimca, G. Sticlaru, Koszul complexes and pole order filtrations, Proc. Edinburg. Math. Soc. 58(2015), 333--354.


\bibitem{DStExpo} A. Dimca, G. Sticlaru, On the exponents of free and nearly free projective plane curves, Rev. Mat. Complut. 30(2017), 259--268.





\bibitem{DStRLM} A. Dimca, G. Sticlaru, Line and rational curve arrangements, and Walther's inequality, arXiv:1803.05386, to appear in Atti Accad. Naz. Lincei Rend. Lincei Mat. Appl. 


\bibitem{DStJump} A. Dimca, G. Sticlaru, On the jumping lines of bundles of logarithmic vector fields along plane curves, arXiv: 1804.06349.

\bibitem{DSt3syz} A. Dimca, G. Sticlaru, Plane curves with three syzygies, minimal Tjurina curves, and nearly cuspidal curves, arXiv: 1810.11766.


\bibitem{DS}
A. Dimca and G. Sticlaru, 
Jacobian syzygies and plane curves with maximal global Tjurina numbers. 
arXiv:1901.05915.  




\bibitem{dPW} A.A. du Plessis,  C.T.C. Wall, Application of the theory of the discriminant to highly singular plane curves, \textit{Math. Proc. Camb.
Phil. Soc.},  \textbf{126} (1999), 259-266. 






\bibitem{E} Ph. Ellia,  
Quasi complete intersections and global Tjurina number of plane curves, J. Pure Appl. Algebra (2019), https://doi.org/10.1016/j.jpaa.2019.05.014.














 
 
 \bibitem{MV19} S. Marchesi, J. Vall\` es, Terao's conjecture for triangular arrangements, 
arXiv:1903.08885.




\bibitem{OT} P. Orlik and H. Terao, {\em Arrangements of Hyperplanes,} Springer-Verlag, Berlin Heidelberg New York, 1992.









\bibitem{STY}  H. Schenck, H. Terao, M. Yoshinaga,  Logarithmic vector fields for curve configurations in $\PP^2$ with quasihomogeneous singularities, Math. Res. Lett. 25 (2018), 1977--1992.






\bibitem{ST} A. Simis, S.O. Toh\u aneanu, Homology of homogeneous divisors, Israel J. Math. 200 (2014), 449-487.



\bibitem{T1}
H. Terao, 
Arrangements of hyperplanes and their freeness I, II. 
\textit{J. Fac. Sci. Univ. Tokyo} \textbf{27} (1980), 293--320.   

\bibitem{W} M. Wakefield,  Derivation degree sequences of non-free arrangements, Communications in Algebra, 47:7(2019), 2654--2666, DOI: 10.1080/00927872.2018.1536207.
	

\bibitem{Y2}
M. Yoshinaga, 
On the freeness of 3-arrangements. 
\textit{Bull. London Math. Soc.} \textbf{37} (2005), no. 1, 126--134. 



\bibitem{Y3}
M. Yoshinaga, 
Freeness of hyperplane arrangements and related topics. 
\textit{Annales de la Faculte des Sciences de Toulouse}, \textbf{23} (2014), no. 2,  483--512. 

\bibitem{Z}
G. M. Ziegler, 
Multiarrangements of hyperplanes and their freeness.  Singularities (Iowa City, IA, 1986),  345--359,


\bibitem{Zi} G. Ziegler, Combinatorial construction of logarithmic differential forms, Adv. Math. 76 (1989), 116-154.

\end{thebibliography}
\end{document}